# DYNAMIC PHYSICAL SYSTEMS : ENERGY BALANCES AND STABILITY ISSUES

M. De la Sen

Department of Electricity and Electronics

Faculty of Sciences. University of Basque Country

Campus de Leioa . Bizkaia. Aptdo. 644 de Bilbao

48080- Bilbao. SPAIN

**Abstract**. This paper relates properties of operators with the well- known concepts of positive realness and passivity properties in dynamic systems and their associate transfer functions. Those concepts together with very close related ones are first examined from a physical point of view. Then , they are related with hyperstability and properties of transfer functions while the hyperstability theorem is revisited and interpreted. Finally, the above concepts are compared to the mathematical concepts of positivity and closely related ones in operator theory in Hilbert spaces.

## I. Introduction

Stability properties of nonlinear dynamic systems have been widely studied in the literature, [1-3], [9-12]. Related properties include, for instance, Lyapunov´s stability /asymptotic stability, absolute stability (i.e. global Lyapunov´s asymptotic stability in the presence of nonlinear static devices belonging to prescribed sectors in the feedback law) or hyperstability / asymptotic hyperstability (i.e. global Lyapunov´s stability / asymptotic hyperstability in the presence of any nonlinear and /or time-varying devices whose time input-output integral satisfy Popov´s type inequalities). While Lyapunov´s stability may be local around the equilibrium, absolute stability/hyperstability are always global in the whole state space and established as a generic property for a set (not just for a single element) of feedback devices for a given forward device or plant. An important physical property is that a positive dynamic system being hyperstable (roughly speaking positive) which is feedback connected with any class of devices satisfying a Popov´s -type inequality implying lower bounding by a negative finite constant is globally Lyapunov´s stable since its input-output energy is nonnegative and bounded for all time, [4-8]. On the other hand, hyperstability for a set of nonlinear/time-varying devices satisfying a certain Popov´s inequality includes the absolute stability of any static nonlinear device that satisfies such an inequality. The above concepts are very related to the more general one of passivity. In an operator theoretical framework, there are well-known related concepts based on positivity of operators, [1]. In this paper, we analyze and inter-relate the various concepts of passivity, hyperstability, positivity, dissipation, conservation, regeneration etc. in Physics from their implications in input-output or power energy balances as well as their strict- type version. We interpret those concepts in a feedback framework related to general stability properties (or roughly speaking hyperstability). Then, we relate those concepts to close properties in the operator theoretical framework formulated in an appropriate Hilbert space.

## II. Physical Concepts Related to Power and Energy Balances

Consider a scalar (only for purposes of facilitating the mathematical treatment and exposition) dynamic systems with instantaneous real input and output signals at time t being, respectively, u (t) and y (t) , then of supplied power (u.(t) y (t) ), whose stored energy and dissipated energy are respectively given by functions S (t) and D(t). Thus, the instantaneous power balance at time $t \geq 0$ and the energy balance in the time interval $[0, t]$ are given , respectively, by :

**Power balance at time t**: $u(t) y(t) = \dot{S}(t) + \dot{D}(t)$

(1.a)

**Energy balance in the time interval $[0, t]$**:

$< u, y >_t = S(t) + D(t) - S(0) - D(0)$ (1.b)

where the dot superscript denotes time-derivative, as usual, $< u, y >_t$ is an abbreviation for a time- integral product (i.e. a scalar product, denoted by $< u, y >_t$, of square-integrable functions u(t) and y (t) on $[0, t]$; i. e. belonging to $L_2[0, t]$) meaning $< u, y >_t = \int_0^t u(\tau) y(\tau) d\tau$. It the time subscript " t " is dropped out from the scalar product definition then the time integral, provided to exist, is extended to infinity; i. e. $< u, y > = \int_0^\infty u(\tau) y(\tau) d\tau$. Note that if truncated input and output signals $u_t$ and $y_t$ replace u and y where $z_t = z(\tau)$ for all $\tau \in [0, t]$ and $z_t = 0$ otherwise in the real axis then $< u, y >_t = < u_t, y_t > = \int_{-\infty}^\infty u_t(\tau) y_t(\tau) d\tau$; i.e. , the input/output energy time-integral may be extended from minus infinity to infinity when using truncated input/output signals. This allows to describe the supplied energy equivalently in the frequency domain via Parseval´s theorem for all finite time even if the input / output product is not potentially square - integrable on [ 0, ∞ ). In the following, we drop the time argument t in order to simplify the notation when no confusion is expected. In the context of dynamic systems, we manipulate a set of energetic-related concepts saying that the system is at time $t \geq 0$ (the constraint $t > 0$ for time is stated explicitly when applicable) :

**a) Regenerative** if it does not dissipate energy but it supplies it to the network. Thus, $\dot{D}(t) < 0$ and $D(t) < D(0)$ so that $u(t) y(t) < \dot{S}(t)$ and $<u, y>_t < S(t) - S(0) < S(t)$. If, in addition, the stored energy decreases with time then $S(t) \leq S(0)$ for all $t \geq 0$ and then $<u, y>_t < 0$ for all $t > 0$.

**b) Passive or Dissipative** if it has energetic losses since $\dot{D}(t) \geq 0$. Thus, $D(t) \geq D(0)$ so that $u(t) y(t) \geq \dot{S}(t)$, and $<u, y>_t \geq S(t) - S(0)$
$$\geq \beta := \min_{t \geq 0} S(t) - S(0) \geq -S(0)$$

Note that β is a real number whose sign depends on each particular situation related to the system's properties. For instance, if $S(t)$ tends asymptotically to zero then $\beta = -S(0)$. However, β is nonnegative ( positive for any $t > 0$ ) if $S(t) \geq S(0)$ ($S(t) > S(0)$ for any $t > 0$). The system is said to be **Strictly Passive or Strictly Dissipative** if $\dot{D}(t) > 0$ for all finite time so that $<u, y>_t > S(t) - S(0)$ for all $t > 0$ except possibly at a set of zero measure. A more complete classification of passivity may be made as follows:

- The system is **Weakly Passive** ( then called Positive as well) if $<u, y>_t \geq 0$ for all $t \geq 0$.

- The system is **Weakly Strictly Passive** (then called Weakly Strictly Positive as well) if $<u, y>_t > 0$ for all $t > 0$.

- The system is **Strongly Strictly Passive** (then called Strongly Strictly Positive as well) if $<u, y>_t > \beta <u, u>_t$ for some real constant $\beta > 0$ and all $t \geq 0$.

**c) Conservative** if $\dot{S}(t) = 0$ ; i.e. the stored energy is kept constant while the supplied energy is entirely dissipated so that :

$$<u, y>_t \geq D(t) - D(0) \geq -D(0).$$

**d) Positive (Strictly Positive)** if $u(t) y(t) \geq 0$ so that $<u, y>_t \geq 0$ ($u(t) y(t) > 0$ and $<u, y>_t > 0$ for all $t > 0$). The specifications **Weakly** or **Strongly** may be used in the same contexts and meanings as for Strict Passivity so that Strictly Positive systems may be specified as Weakly Strictly Positive or Strongly Strictly Positive ones, respectively. Positive systems may be equivalently named as Weakly Passive Systems.

**e) It satisfies Popov's Inequality**. If for some finite real constant $\gamma_0$ and all $t \geq 0 <u, y>_t \geq -\gamma_0^2 > -\infty$.

**Remarks : (1)** The above concepts may also be applicable only to some finite time subinterval $[t_1, t_2]$ in such a way that the system may be characterized under different properties in the above context through time.

**(2)** Both Passive and Positive dynamic Systems satisfy Popov's Inequality.

**(3)** A system which satisfies Popov's Inequality is always passive or conservative but not necessarily Positive (i.e., not necessarily Weakly Passive).

**(4)** If a system is regenerative and $S(t) \leq S(0)$, for all finite time, the energy supplied is negative for all finite time so that in fact the system supplies energy to the connected network. Also, its supplied input/output energy is upper-bounded by a negative real number.

**(5)** A system is both Passive and Positive if $<u, y>_t \geq \beta \geq 0$. A system is Passive but not Positive( then not Weakly Passive) in some interval [ 0, t] if there exists a finite negative β such that $<u, y>_t \geq \beta$. Then, the system satisfies Popov's Inequality as well.

### III. Hyperstability

The above concepts play a crucial role in the properties of hyperstability and asymptotic hyperstability which, as stated in the introduction, generalize the concept of absolute stability which, on the other hand, generalizes the standard one of global Lyapunov's stability. Assume a negative feedback configuration where the forward loop is defined by a linear time-invariant input/output operator (or plant) from the input space to the output space $G : U \to Y$ while the feedback loop is a, in general, nonlinear and / or time-varying operator (or feedback controller) whose output space is equal to the input space to the forward loop $F : Y \to V \equiv U$ such that if u is in U then v = -u is in V identical to U. **Assume that the G-operator is Strictly Positive and the feedback one is anyone satisfying a Popov's -type Inequality so that** :

$$<u, y>_t \geq 0 ; -<u, y>_t = <v, y>_t \geq -\gamma_0^2 > -\infty$$
(2)

Combining the above two relationships, one gets that the supplied input/output energy during the time interval [0, t] satisfies after using Parseval's theorem and assuming that the input is not identically zero within such an interval:

$$E(t) = <u, y>_t = <u_t, y_t>$$
$$= <u_t, g * u_t> = <u_t, h u_t>$$
$$= (2\pi)^{-1} <\hat{u}_t, \hat{y}_t>$$
$$= (2\pi)^{-1} <\hat{u}_t, \hat{g} \hat{u}_t> \qquad (3)$$

where j is the imaginary unit, the symbol * denotes the convolution integral, g and $\hat{g}$ being the impulse response and the frequency response ( i.e. its Fourier transform $F(.)$) associated with the physical filter of the forward input-output G-operator, and h being a time operator from U to Y defining the convolution integral in the time-domain, namely:

$$g * u_t = h(u_t)(t) = \int_{-\infty}^{\infty} g(\tau) u_t(t-\tau) d\tau$$

$$= \int_0^t g(\tau) u(t-\tau) d\tau$$

$$\hat{u}_t(j\omega) = F(u_t) = \int_{-\infty}^{\infty} u_t(\tau) e^{-j\omega\tau} d\tau$$

such Fourier transforms always exist for finite time since the corresponding integrals exist. Note that the input/output energy is expressed equivalently in the time-domain (first line of identities in eqn. 3) and in the frequency domain (second line of identities in eqns. 3). Thus,

$$E(t) = (2\pi)^{-1} \int_{-\infty}^{\infty} \hat{u}_t(j\omega)(\hat{g}(j\omega)\hat{u}_t(-j\omega)) d\omega$$

$$= (2\pi)^{-1} \int_{-\infty}^{\infty} \operatorname{Re} \hat{g}(j\omega) |\hat{u}_t(j\omega)|^2 d\omega$$

$$= (2\pi)^{-1} <\hat{u}_t, (\operatorname{Re} \hat{g})\hat{u}_t> \quad (4)$$

with the last inner product being defined in the frequency input / output spaces by using the identities (3) where the odd symmetry property of the imaginary part of the hodograph $\operatorname{Im}(\hat{g}(j\omega)) = -\operatorname{Im}(\hat{g}(-j\omega))$ has been used.

### A) Asymptotic Hyperstability for Strongly Strictly Positive Real transfer functions

Now, if the h and Re $\hat{g}$ are Strictly Positive ( or, in particular, **Strongly Strictly Passive**) operators then $d = \operatorname*{Min}_{\omega \geq 0} \operatorname{Re} \hat{g}(j\omega) > 0$ [checking for negative frequencies is not necessary since $\operatorname{Re}(\hat{g}(j\omega)) = \operatorname{Re}(\hat{g}(-j\omega))$]. It is then said that **the transfer function $\hat{g}(s)$ is Strongly Strictly Positive Real**, i.e. $\operatorname{Re} \hat{g}(s) > d \geq 0$ for $\operatorname{Re} s \geq 0$ so that $\operatorname{Re} \hat{g}(j\omega) \geq d > 0$ for all real $\omega$ [4-8], so that one gets directly from (4) combined with the second relationship in (2) for the feedback loop :

$$\infty > \gamma_0^2 \geq E(t)$$
$$\geq (2\pi)^{-1} d \int_{-\infty}^{\infty} |\hat{u}_t(j\omega)|^2 d\omega = d \int_0^t u^2(\tau) d\tau > 0$$
for $t > 0$ \quad (5)

so that taking limits as $t \to \infty$ it follows that the input is bounded for all time and it converges to zero asymptotically continuous (or it only has bounded isolated discontinuities). Since $\hat{g}(s)$ is Strongly Strictly Positive Real then it is strictly stable (i.e. its poles have negative real parts) and non-strictly proper (i.e. it has the same number of poles and zeros -or relative degree zero). Its inverse $1/\hat{g}(s)$ is also Strongly Strictly Positive Real, strictly stable and non-strictly proper but proper (and then realizable) so that $1/d = \operatorname*{Min}_{\omega \geq 0}(\hat{g}^{-1}(j\omega)) > 0$. Thus, (5) might be re-arranged by using $\hat{u}(j\omega) = \hat{g}^{-1}(j\omega)\hat{y}(j\omega)$ as follows for $t > 0$ :

$$\infty > \gamma_0^2 \geq E(t)$$
$$\geq (2\pi)^{-1} d^{-1} \int_{-\infty}^{\infty} |\hat{y}_t(j\omega)|^2 d\omega = d \int_0^t y^2(\tau) d\tau > 0$$
\quad (6)

Then, taking limits as above as time tends to infinity, one concludes that the output is bounded provided that it is continuous almost everywhere and tends asymptotically to zero. The asymptotic hyperstability theorem is formulated as follows, [5] .**Thus, if the plant is Strongly Strictly Passive (so that its transfer function is Strongly Strictly Positive Real) while the feedback loop is anyone satisfying a Popov´s type Inequality then the closed-loop system is asymptotically hyperstable ( i.e. globally Lyapunov´s asymptotically stable for the class of feedback laws satisfying the Popov´s Inequality in (2).** If the transfer function is Weakly Strictly Positive Real, so that its associate time and frequency domain operators are Weakly Strictly Passive, then $\operatorname{Re} \hat{g}(j\omega) > 0$ for all finite $\omega$ but $\lim_{\omega \to \pm\infty} \operatorname{Re} \hat{g}(j\omega) = 0$.

### B) Asymptotic Hyperstability for Weakly Strictly Positive Real transfer functions

Thus, the above reasoning needs to be modified to get the asymptotic hyperstability result. **Assume that the transfer function is Weakly Strictly Positive Real** with $\operatorname{Re} \hat{g}(j\omega) > 0$ for all finite $\omega$ and $\lim_{\omega \to \pm\infty} \omega^2 \operatorname{Re} \hat{g}(j\omega) \geq d_0 > 0$. Then, we perform multiplication and division by the squared-frequency in the frequency domain integrals of (5) to get instead:

$$\infty > \gamma_0^2 \geq E(t)$$
$$\geq (2\pi)^{-1} d_0 \int_{-\infty}^{\infty} |\hat{\delta}_t(j\omega)|^2 d\omega$$
$$= (2\pi)^{-1} d_0 \int_{-\infty}^{\infty} \left(\frac{\hat{u}_t(j\omega)}{j\omega}\right)\left(\frac{\hat{u}_t(-j\omega)}{-j\omega}\right) d\omega$$
$$= d_0 \int_0^t \delta^2(\tau) d\tau > 0 \text{ for } t>0 \quad (7)$$

where $\delta(.)$ is the input time-integral. Thus, it follows that this integral converges to zero as time tends to infinity so that the input should exhibit that limit behavior. Continuing with such a development one gets the following conclusion. **Thus, if the plant is Weakly Strictly Passive (so that its transfer function is Weakly Strictly Positive Real) while the feedback loop is anyone satisfying a Popov´s type Inequality then the closed-loop system is asymptotically hyperstable (i.e. globally Lyapunov´s asymptotically stable for the class of feedback laws satisfying the Popov´s Inequality in (2).**

### C) Further Comments

Note also that, in both cases of Strict Positive Realness , the plant input/ output energy and supplied power are at the same time positive and bounded for all time : i.e, bounded above with a finite bound and strictly positively bounded from below for all time**.**

A key associate property is that the absolute maximum input/output phase deviation is 90º and that the system is strictly stable of strictly stable inverse in the case of strict positivity or passivity and critically stable (of inverse being critically stable as well) with eventual simple imaginary poles of nonnegative associate residuals. Also, the hodographs of frequency responses are confined to the first and third quadrants of the complex plane and they are never tangent to the imaginary axis if the system is Strongly Strictly Positive Real. Note that another important aspect is the role played by the feedback device. Note that while the forward loop is strictly positive / passive (and then dissipative) the feedback one might have negative supplied energy (at least during certain time intervals) so that it may be regenerative at least during certain time intervals. In this case, the upper-bound of the feedback input/output integral satisfying Popov´s Inequality is a negative real number during such time intervals. This leads to the weaker sufficient conditions for achieving closed-loop stability, when adopting a physical point of view concerning weakness of dynamics constraints, but it is not always the case concerning the fulfillment of Popov´s Inequality. For instance, if the feedback loop consists of a dynamics-free nonlinearity inside the first/third quadrants, as in the standard absolute stability problem, then the above mentioned upper-bound is always positive for the scalar product satisfying a Popov´s Inequality type lower-bound what means that the feedback device is either conservative or dissipative as it is the forward device (plant) while maintaining closed-loop stability in terms of hyperstability.

We can also point out by using again Parseval´s theorem in (4) to interpret it in the time-domain via the bounds in (5) that $\infty > \int_0^t g(\tau) u^2(\tau) d\tau > 0$ if the system is Strictly Passive (or Strictly Positive), so that its transfer function is Strictly Positive Real and $\infty > \int_0^t g(\tau) u^2(\tau) d\tau \geq 0$ for all $t > 0$ if the system is Weakly Passive (or Positive), so that its transfer function is Positive Real. As a result, the impulse response g (t) is a strictly positive function and bounded above for all time t > 0 if the system is either Weakly or Strongly Strictly Passive/ Positive and g (t) ≥ 0 and bounded above for all time t > 0 if the system is Weakly Passive/ Positive. If the system is only Positive/ Weakly Passive then g(t) does not converge asymptotically to zero. Thus, the last inequality ensures that the input u (t) is bounded. Since the transfer function is (perhaps critically) stable [since Positive Real] then the output is bounded as well and (in general, non asymptotic) hyperstability is guaranteed.

**D) Asymptotic Hyperstability for Positive Real transfer function with a single Pole at the Origin (Popov´s Simplest Particular Case)**

Now, assume the case that the plant input is not trivially zero and the forward loop is only (nonstrict) Positive /Weakly Passive while its transfer function possess only a single pole at s = 0. Assume also that $\hat{g}_1(s) = s\hat{g}(s)$ is Strictly Positive Real. After relating real and imaginary parts of $\hat{g}(s)$ and $\hat{g}_1(S)$, one gets

$$\text{Re } \hat{g}(j\omega) = \frac{\text{Im } \hat{g}_1(j\omega)}{\omega}$$

and Re $\hat{g}_1(j\omega) = -\omega \text{Im}\hat{g}(j\omega)$ so

Im $\hat{g}(j\omega) \leq 0$ and Im $\hat{g}_1(j\omega) \leq 0$ for $\omega \geq 0$ should hold in addition. Now, note that

$$E(t) = (2\pi)^{-1} \int_{-\infty}^{\infty} \text{Re } \hat{g}_1(-j\omega) \hat{\bar{u}}_t(j\omega) d\omega$$

where $d_1 > 0$ provided that $\hat{g}_1(s) = s\hat{g}(s)$ is Strongly Strictly Positive Real (so that strictly stable and of relative degree zero or plus unity) since $\hat{g}(s)$ is Positive Real with a single pole at s=0 and

$$\hat{\bar{u}}_t(j\omega) = \frac{|\hat{u}_t(j\omega)|^2}{j\omega} \text{ so that}$$

$$\hat{\bar{u}}_t(\tau) = \int_{-\infty}^{\infty} \frac{|\hat{u}_t(j\omega)|^2}{j\omega} e^{j\omega\tau} d\omega$$

$$E(t) \geq (2\pi)^{-1} d_1 \int_{-\infty}^{\infty} |\hat{\delta}_t(j\omega)| |\hat{u}_t(-j\omega)| d\omega$$

$$= d_1 \int_{-\infty}^{\infty} \delta_t(\tau) |u_t(\tau)| d\tau = d_1 \int_0^t \delta(\tau) |u(\tau)| d\tau$$
$$> 0$$

for any nontrivial input where $\delta(t) = \int_0^t |u(\tau)| d\tau$.

After combining the above inequalities with Popov´s Inequality of the feedback device, one gets that the input is bounded, square-integrable and converges to zero. The output has the same properties since $\hat{g}_1$ is strongly positive real. Then, asymptotic hyperstability follows also in this particular case of (nonstrict) positive realness where $\hat{g}_1(S)$ is Strongly Strictly Positive Real. The proof for the case when $\hat{g}_1(S)$ is Weakly Strictly Positive Real is quite similar but more involved and it may be addressed by proceeding with $\hat{g}_1(S)$ as in the case of Weakly Strictly Positive Real transfer function discussed previously in the context of asymptotic hyperstability for strict realness of the forward loop. A very related case is that the Simplest Particular Case (i.e. Positive Realness of the plant with a single pole at the origin) leads to absolute stability (global asymptotic Lyapunov´s stability) for any nonlinear device which only generates a zero output when its input takes a zero value.

### IV. Links with Operator Theory

All the above results may be interpreted in the context of operators. We consider the input and Output spaces U ( identical to V ) and Y as Hilbert linear subspaces ( i.e. Banach spaces , namely, normed spaces where any Cauchy sequence has a limit in those spaces) of the set or real square -integrable functions $L_2 \equiv L_2(0, \infty)$ endowed with the inner product (semi) norm ; i.e. if $u \in U$ then $|u| = \sqrt{<u, u>}$ and a similar norm is defined for the output signal on Y . Since , we have to deal with limits as time tends to infinity, it cannot be " a priori" guaranteed that the input/ output functions are square-integrable over $(0, \infty)$ since this has been a previous issue in the stability proofs of the former section. Therefore, the formalism is more properly established on
$$L_{2e} := \{ f : [0, \infty) \to R \ / \ f_t \in L_2 \ \forall t \in [0, \infty) \}$$
$$\equiv \bigcup_{0 \leq t < \infty} (L_2[0, t])$$

i.e., the set of square-integrable truncated functions for any finite truncation time. Thus, for all finite time, we can consider the (truncated) input and output signals of the dynamic system as members of that set. Also, since the $L_2$-norm is rather a seminorm , since it is defined through an integral, we consider as identical all input and output signals belonging to classes that only differ possibly on sets of zero measure of $(0, \infty)$ . Now, we pay our attention to a key identity recovered from (3), namely,

$E(t) = \ <u_t, h u_t> \ \geq 0$ for all $t \geq 0$

for all finite t . In our context, we say that this holds for any $u_t \in L_2$ for finite time (which, in fact, is identical to say for any $u \in L_2[0, t]$ for any finite time). That means that the Convolution Operator is Positive if the transfer function of the plant is Positive Real or Strictly Positive Real. That leads, trough Parseval´s theorem, to the fact that the associate response frequency operator which is the mapping between the corresponding input and output frequency linear spaces ( being identified in particular with the real part of the frequency response $\text{Re} \hat{g}(j\omega)$)) is also positive, respectively, strictly positive. Positive Operators are self-adjoin operators. If the two-sided boundedness of the input/output energy balance discussed in the above section (finite above and below strictly from zero) holds for all time, which requires for the feedback loop to satisfy Popov´s Inequality, then the system is asymptotically hyperstable since we can take limits as time tends to infinity to conclude that $u \in L_2$ , u tends to zero as time tends to infinity while it is bounded for all time , provided that $u_t \in L_2$. In order to interpret all the results of the previous sections in the context of operator theory, we can extend the definition of positive operators to passive ones together with their strict versions as follows:

The h-operator is (and so it is the operator $\hat{h}(j\omega) := \text{Re} \ \hat{g}(j\omega)$ through Parseval´s theorem) :

. **Passive or Dissipative** : $<u_t, h u_t> \ \geq \beta$ for some real constant $\beta$ all $t \geq 0$. This implies $<\hat{u}_t, \hat{h} \hat{u}_t> \ \geq \beta$.

. **Positive** if $\beta = 0$; **Weakly Strictly Positive / Passive** if $<\hat{u}_t, \hat{h} \hat{u}_t> \ \geq \beta > 0$ for all nonzero $u_t$ and all $t > 0$; and **Strongly Strictly Positive / Passive if** $<u_t, h u_t> \ \geq \beta <u_t, u_t>$ with $\beta > 0$ for all $t \geq 0$. Since the properties of the h-operator induce similar properties on the $\hat{h}$-operator, it follows that :

. If h is Positive then $\hat{h}$ is positive as well, $\text{Re} \ \hat{g} \geq 0$ for all real ω so that $\hat{g}$ is Positive Real. As a result, it is (perhaps critically) stable with relative degree zero or plus unity ( if realizable) with residuals at the critically stable (necessarily simple) poles ( if any) being nonnegative, having inverse Positive Real and producing an absolute input / output phase deviation of at most 90º.

. If h is Strongly Strictly Positive then $\hat{h}$ is Strongly Strictly Positive as well, $\text{Re} \ \hat{g} > 0$ for all real ω so that $\hat{g}$ is Strongly Strictly Positive Real. As a result, it is strictly stable with relative degree zero, having inverse Strictly Positive Real and producing an absolute input / output phase deviation of at most 90º.

. If h is weakly Strictly Positive then $\hat{h}$ is weakly Strictly Positive as well, $\text{Re} \ \hat{g} > 0$ for all real finite ω [with $\hat{g}$ tending to zero as the absolute frequency tends to infinity and $\omega^2 \hat{g}$ tending to a positive number as the absolute frequency tends to infinity] so that $\hat{g}$ is weakly Strictly Positive Real. As a result, it is strictly stable, having inverse Strictly Positive Real and producing an absolute input / output phase deviation of at most 90º.The proof of asymptotic hyperstability requires that the feedback F- operator satisfy Popov´s Inequality and such a proof is addressed as indicated in the previous section.

ACKNOWLEDGMENTS
The author is very grateful to MCYT by its partial support of this work through Project DPI 2003-0164 and to UPV/EHU by its support through Grant 9 / UPV 00I06.I06-15263/2003.